\documentclass[12pt]{amsart}
\usepackage{amsfonts}
\usepackage{amsmath,amssymb,amsthm}

\newcommand{\R}{\mathbb R}

\renewcommand{\span}{\mathrm{span}}
\newcommand{\tr}{\mathrm{tr}}

\newtheorem{thm}{Theorem}[section]

\newtheorem{prop}[thm]{Proposition}
\theoremstyle{definition}
\newtheorem{defn}[thm]{Definition}
\theoremstyle{remark}

\newcommand{\ds}{\displaystyle}

\newcommand\topstrut{\rule{0mm}{2.9ex}}
\newcommand\bottomstrut{\rule[-1.5ex]{0mm}{1.5ex}}
\newcommand\titlestrut{\topstrut\bottomstrut}

\begin{document}

\title[TIMELIKE SURFACES WITH ZERO MEAN CURVATURE IN  MINKOWSKI 4-SPACE]
{TIMELIKE SURFACES WITH ZERO MEAN CURVATURE IN  MINKOWSKI 4-SPACE}

\author{Georgi Ganchev and Velichka Milousheva}
\address{Bulgarian Academy of Sciences, Institute of Mathematics and Informatics,
Acad. G. Bonchev Str. bl. 8, 1113 Sofia, Bulgaria}
\email{ganchev@math.bas.bg}

\address{Bulgarian Academy of Sciences, Institute of Mathematics and Informatics,
Acad. G. Bonchev Str. bl. 8, 1113, Sofia, Bulgaria; "L. Karavelov"
Civil Engineering Higher School, 175 Suhodolska Str., 1373 Sofia,
Bulgaria}
\email{vmil@math.bas.bg}

\subjclass[2000]{Primary 53A07, Secondary 53A10}

\keywords{timelike surfaces with zero mean curvature, canonical parameters, natural
partial differential equations, rotational surfaces of Moore type}

\begin{abstract}
On any  timelike surface with zero mean curvature in the four-dimensional Minkowski
space we introduce special geometric (canonical) parameters and prove that the
Gauss curvature and the normal curvature of the surface satisfy a system of two natural
partial differential equations. Conversely, any two solutions to this system
determine a unique (up to a motion) timelike surface with zero mean curvature  so that the given parameters
are canonical. We find all timelike surfaces with zero mean curvature  in the
class of rotational surfaces of Moore type. These examples give rise to a one-parameter
family of solutions to the system of natural partial differential equations
describing  timelike surfaces with zero mean curvature.
\end{abstract}

\maketitle

\section{Introduction} \label{S:Intr}

Studying minimal surfaces in the four-dimensional Euclidean space $\R^4$, Itoh proved
in \cite{Itoh} that any minimal non-superconformal surface $M^2$ admits locally special
isothermal parameters. On the base of this result, de Azevero Tribuzy and Guadalupe
proved in \cite{Guad-Trib} the following:
\vskip 1mm
{\it The Gauss curvature $K$ and the curvature of the normal connection $\varkappa$
(the normal curvature) of a minimal non-superconformal surface, parameterized by
special isothermal parameters, satisfy the following system of  partial differential
equations:
\begin{equation}\label{Eq-0}
\begin{array}{l}
\ds{(K^2 - \varkappa^2)^{\frac{1}{4}}\, \Delta \ln |\varkappa - K|} = 2(2K - \varkappa);\\
[2mm]
\ds{(K^2 - \varkappa^2)^{\frac{1}{4}}\, \Delta \ln |\varkappa + K| = 2(2K + \varkappa)}.
\end{array}
\end{equation} Conversely, any solution ($K$, $\varkappa$) to system \eqref{Eq-0} determines uniquely (up to a
motion in $\R^4$) a minimal non-superconformal surface with Gauss curvature $K$
and normal curvature $\varkappa$.}

\vskip 2mm
In the present paper we study the local theory of timelike surfaces with zero mean curvature in the four-dimensional Minkowski space $\R^4_1$.

In Section \ref{S:timelike}  we introduce (Theorem \ref{T:Canonocal parameters-a}) special isothermal parameters (canonical parameters)
on any  timelike surface with zero mean curvature  free of flat points and clear the geometric nature
of these parameters.

Using Theorem \ref{T:Canonocal parameters-a}, we prove the fundamental Theorem \ref{T:Fundamental Theorem 2 - timelike},
which  states as follows:
\vskip 1mm
{\it Let $M^2$ be a  timelike surface  with zero mean curvature free of flat points.
If $M^2$ is parameterized by canonical parameters, then the Gauss curvature $K$ and the normal curvature $\varkappa$
satisfy the system of natural partial differential equations:
\begin{equation}\notag
\begin{array}{l}
\ds{(K^2 + \varkappa^2)^{\frac{1}{4}}\, \Delta^h \ln (K^2 + \varkappa^2)^{\frac{1}{8}}} = K ;\\
[2mm]
\ds{(K^2 + \varkappa^2)^{\frac{1}{4}}\, \Delta^h \arctan \frac{\varkappa}{K} = 2 \varkappa}.
\end{array}
\end{equation}

Conversely, any solution ($K$, $\varkappa$) to the above system,
determines a unique (up to a motion in $\R^4_1$)  timelike  surface  with zero
mean curvature free of flat points  with Gauss curvature $K$
and normal curvature $\varkappa$. Moreover the given parameters
are canonical for this surface.}

In order to apply the above results, we consider the class of timelike rotational
 surfaces of Moore type.
In Section \ref{S:Examples}  we find the meridian curves of all  timelike
rotational surfaces of Moore type with zero mean curvature.
We note that for this class of surfaces the functions $K$ and $\varkappa$  are only functions of one variable.

In the last section we make a survey of
the background systems of partial differential equations describing  surfaces with zero mean curvature
in $\R^4$ or  $\R^4_1$. These systems written in canonical form are given below:

\vskip 3mm
\centerline{\emph{Surfaces with zero mean curvature}:}

\vskip 2mm
\begin{center}
\begin{tabular}{|c|c|c|}
\hline
\titlestrut
\quad \emph{surfaces in} $\R^4$
 \quad & \quad \emph{spacelike surfaces in} $\R^4_1$ \quad &  \quad \emph{timelike surfaces
in} $\R^4_1$ \quad \\

\hline
\titlestrut

$\begin{array}{l}
\\
\Delta X = 2 e^{X} \cosh Y\\
[2mm]
 \Delta Y = 2 e^{X} \sinh Y
 \\[3mm]
\end{array}$

& $\begin{array}{l}
\\
\Delta X = 2 e^{X} \cos Y\\
[2mm]
\Delta Y = 2 e^{X} \sin Y
\\[3mm]
\end{array}$

 & $\begin{array}{l}
 \\
\Delta^h X = 2 e^{X} \cos Y\\
[2mm]
\Delta^h Y = 2 e^{X} \sin Y
\\[3mm]
\end{array}$\\

\hline

\end{tabular}
\end{center}

\section{Preliminaries} \label{S:Pre}

Let  $\R^4_1$ be the four-dimensional  Minkowski space endowed with the metric
$\langle , \rangle$ of signature $(3,1)$ and $Oe_1e_2e_3e_4$ be a
fixed orthonormal coordinate system in $\R^4_1$, i.e. $e_1^2 =
e_2^2 = e_3^2 = 1, \, e_4^2 = -1$, giving the orientation of
$\R^4_1$. The standard flat metric is given in local coordinates by
$dx_1^2 + dx_2^2 + dx_3^2 -dx_4^2.$

A surface $M^2$ in $\R^4_1$ is said to be
\emph{spacelike} if $\langle , \rangle$ induces  a Riemannian
metric $g$ on $M^2$. Thus at each point $p$ of a spacelike surface
$M^2$ we have the following decomposition
$$\R^4_1 = T_pM^2 \oplus N_pM^2$$
with the property that the restriction of the metric
$\langle , \rangle$ onto the tangent space $T_pM^2$ is of
signature $(2,0)$, and the restriction of the metric $\langle ,
\rangle$ onto the normal space $N_pM^2$ is of signature $(1,1)$.

A surface $M^2$ in $\R^4_1$ is said to be
\emph{timelike} if the induced metric $g$ on $M^2$ is a metric with index 1, i.e.
at each point $p$ of a timelike surface $M^2$ we have the following decomposition
$$\R^4_1 = T_pM^2 \oplus N_pM^2$$
with the property that the restriction of the metric
$\langle , \rangle$ onto the tangent space $T_pM^2$ is of
signature $(1,1)$, and the restriction of the metric $\langle ,
\rangle$ onto the normal space $N_pM^2$ is of signature $(2,0)$.

Denote by $\nabla'$ and $\nabla$ the Levi Civita connections on $\R^4_1$ and $M^2$, respectively.
Let $x$ and $y$ denote vector fields tangent to $M$ and let $\xi$ be a normal vector field.
Then the formulas of Gauss and Weingarten give  decompositions of the vector fields $\nabla'_xy$ and
$\nabla'_x \xi$ into a tangent and a normal component:
$$\begin{array}{l}
\vspace{2mm}
\nabla'_xy = \nabla_xy + \sigma(x,y);\\
\vspace{2mm}
\nabla'_x \xi = - A_{\xi} x + D_x \xi,
\end{array}$$
which define the second fundamental tensor $\sigma$, the normal connection $D$ and the shape operator $A_{\xi}$ with respect to $\xi$.
The mean curvature vector  field $H$ of the surface $M^2$ is defined as $H = \ds{\frac{1}{2}\,  \tr\, \sigma}$.
Thus, if $M^2$ is a spacelike surface and $\{x,y\}$ is a local orthonormal frame  of the tangent bundle,
the mean curvature vector  field is
$H = \ds{\frac{1}{2} \left(\sigma(x,x) +  \sigma(y,y)\right)}$. If $M^2$ is timelike, then
$H = \ds{\frac{1}{2} \left(- \sigma(x,x) + \sigma(y,y)\right)}$, where   $\{x,y\}$
is a  local orthonormal frame  of the tangent bundle such that $\langle x, x \rangle = - 1$, $\langle y, y \rangle = 1$.

We shall study  timelike surfaces with zero mean curvature, i.e. $H=0$ at each point.

\section{Fundamental theorem for  timelike surfaces with zero mean curvature}\label{S:timelike}

In the local theory of surfaces in Euclidean space  a statement of
significant importance is a theorem of Bonnet-type giving the
natural conditions under which the surface is determined up to a
motion.
In this section we shall introduce canonical parameters  on each  timelike surface free of flat points with $H=0$
and using the canonical parameters we shall  prove a fundamental theorem of Bonnet-type for the class of timelike surfaces with zero mean curvature  free of flat points.

\vskip 2mm
Let $M^2: z=z(u,v), \,\, (u,v) \in \mathcal{D}$ $(\mathcal{D} \subset \R^2)$  be a local parametrization on a
timelike surface in $\R^4_1$.
The tangent space at an arbitrary point $p=z(u,v)$ of $M^2$ is $T_pM^2 = \span \{z_u,z_v\}$. We assume that
$\langle z_u,z_u \rangle < 0$, $\langle z_v,z_v \rangle > 0$.
Hence,
$E(u,v)<0, \; G(u,v)>0$ and we set $W=\sqrt{- EG+F^2}$.
We choose an orthonormal frame field $\{e_1, e_2\}$ of the normal bundle, i.e. $\langle
e_1, e_1 \rangle =1$, $\langle e_2, e_2 \rangle = 1$, $\langle e_1, e_2 \rangle = 0$.
Then we have the following derivative formulas:
$$\begin{array}{l}
\vspace{2mm} \nabla'_{z_u}z_u=z_{uu} = - \Gamma_{11}^1 \, z_u +
\Gamma_{11}^2 \, z_v + c_{11}^1\, e_1 + c_{11}^2\, e_2;\\
\vspace{2mm} \nabla'_{z_u}z_v=z_{uv} = - \Gamma_{12}^1 \, z_u +
\Gamma_{12}^2 \, z_v + c_{12}^1\, e_1 + c_{12}^2\, e_2;\\
\vspace{2mm} \nabla'_{z_v}z_v=z_{vv} = - \Gamma_{22}^1 \, z_u +
\Gamma_{22}^2 \, z_v + c_{22}^1\, e_1 + c_{22}^2\, e_2,\\
\end{array}$$
where $\Gamma_{ij}^k$ are the Christoffel's symbols and the functions $c_{ij}^k, \,\, i,j,k = 1,2$  are given by
$$\begin{array}{ll}
\vspace{2mm}
c_{11}^1 = \langle z_{uu}, e_1 \rangle; & \qquad c_{11}^2 = \langle z_{uu}, e_2 \rangle;\\
\vspace{2mm}
c_{12}^1 = \langle z_{uv}, e_1 \rangle; & \qquad c_{12}^2 = \langle z_{uv}, e_2 \rangle;\\
\vspace{2mm} c_{22}^1 = \langle z_{vv}, e_1 \rangle; & \qquad
c_{22}^2 = \langle z_{vv}, e_2 \rangle.
\end{array} $$

Obviously, the surface $M^2$ lies in a 2-plane if and only if
$M^2$ is totally geodesic, i.e. $c_{ij}^k=0, \; i,j,k = 1, 2.$ So,
we assume that at least one of the coefficients $c_{ij}^k$ is not
zero.

Let us consider the following determinants
$$\Delta_1 = \left|%
\begin{array}{cc}
\vspace{2mm}
  c_{11}^1 & c_{12}^1 \\
  c_{11}^2 & c_{12}^2 \\
\end{array}%
\right|, \quad
\Delta_2 = \left|%
\begin{array}{cc}
\vspace{2mm}
  c_{11}^1 & c_{22}^1 \\
  c_{11}^2 & c_{22}^2 \\
\end{array}%
\right|, \quad
\Delta_3 = \left|%
\begin{array}{cc}
\vspace{2mm}
  c_{12}^1 & c_{22}^1 \\
  c_{12}^2 & c_{22}^2 \\
\end{array}%
\right|.$$
The condition $\Delta_1 = \Delta_2 = \Delta_3 = 0$  characterizes points at which
the space  $\{\sigma(x,y):  x, y \in T_pM^2\}$ is one-dimensional.
In \cite{GM1}, \cite{GM4},  and \cite{GM5} we called such points  \emph{flat points} of the surface.
These points are analogous to flat points in the theory of surfaces in $\R^3$.
In \cite{Lane} and \cite{Little} such points are called inflection points.
The notion of an inflection point is introduced for 2-dimensional surfaces in a 4-dimensional affine space $\mathbb{A}^4$.
E. Lane \cite{Lane} has shown that every point of a surface in $\mathbb{A}^4$ is an inflection point
if and only if the surface is developable or lies in a 3-dimensional space.
Further we consider timelike surfaces free of flat points, i.e. $(\Delta_1, \Delta_2, \Delta_3) \neq (0,0,0)$.

\vskip 2mm
Let  $M^2: z=z(u,v), \,\, (u,v) \in \mathcal{D}$  be a timelike surface free of flat points with zero mean curvature, i.e.
$H = 0$ at each point.
Without loss of generality we assume that the parametrization is orthogonal ($F = 0$) and denote the
unit vector fields $\displaystyle{x=\frac{z_u}{\sqrt{- E}}, \;
y=\frac{z_v}{\sqrt G}}$. The mean curvature vector field is
$H = \ds{\frac{1}{2} \left(- \sigma(x,x) + \sigma(y,y)\right)}$. Since $H = 0$
then $\sigma(x,x) = \sigma(y,y)$ and we have the following decompositions:
\begin{equation}\label{Eq-1-a} \notag
\begin{array}{l}
\vspace{2mm}
\sigma(x,x) = a \,e_1  + b \,e_2; \\
\vspace{2mm}
\sigma(x,y) = c \,e_1 + d \,e_2;  \\
\vspace{2mm}
 \sigma(y,y) =  a \,e_1 + b \,e_2;
\end{array}
\end{equation}
where $a$, $b$, $c$, $d$ are functions on $M^2$.
The assumption that $M^2$ is free of flat points implies $a^2+ b^2 \neq 0$, \,$c^2 +d^2 \neq 0$,\,
$a\,d - b\,c \neq 0$.

Let $\{\overline{x},\overline{y}\}$ be  another orthonormal tangent frame field, given by
$$\begin{array}{l}
\overline{x} = \cosh \varphi \, x + \sinh \varphi \, y,\\
[2mm]
\overline{y} = \sinh \varphi \, x +  \cosh \varphi \, y.
\end{array}$$
Then the corresponding functions
$\overline{a}, \, \overline{b}, \, \overline{c}, \, \overline{d}$ are expressed as follows:
\begin{equation}\label{Eq-2-a}
\begin{array}{ll}
\overline{a} = a \cosh 2\varphi + c \sinh 2\varphi, & \qquad
\overline{b} = b \cosh 2\varphi + d \sinh 2\varphi,\\
[2mm]
\overline{c} =  a \sinh 2\varphi + c \cosh 2\varphi, & \qquad
\overline{d} =  b \sinh 2\varphi + d \cosh 2\varphi.
\end{array}
\end{equation}

We shall find a local orthonormal tangent frame field $\{\overline{x}, \overline{y}\}$, defined in
${\mathcal D}_0 \subset {\mathcal D}$, such that
$\langle \sigma(\overline{x},\overline{x}), \sigma(\overline{x},\overline{y}) \rangle = 0$, or equivalently
$\overline{a}\,\overline{c} + \overline{b}\,\overline{d} = 0$.

Equalities \eqref{Eq-2-a} imply
\begin{equation} \label{Eq-14}
\overline{a}\, \overline{c } + \overline{b}\, \overline{d} = (a\,c + b\,d)\,
\cosh 4 \varphi + \ds{\frac{ a^2 + b^2 + c^2 + d^2}{2}\, \sinh 4 \varphi}.
\end{equation}

If $ a\,c + b\,d = 0$ for all $(u,v) \in {\mathcal D}$ then $\langle \sigma(x,x), \sigma(x,y)\rangle = 0$ at each point of the surface.
If $a\,c + b\,d \neq 0$ at the point $(u_0, v_0) \in {\mathcal D}$, then there exists
${\mathcal D}_0 \subset {\mathcal D}$ ($(u_0, v_0) \in {\mathcal D}_0$), such that
$a\,c + b\,d \neq 0$ for all $(u, v) \in {\mathcal D}_0$.
We denote $A= \ds{- \frac{2(a\,c + b\,d)}{a^2 + b^2 + c^2 + d^2}}$. Then
$A - 1= \ds{- \frac{(a + c)^2 + (b + d)^2}{a^2 + b^2 + c^2 + d^2}}$;
$A + 1= \ds{ \frac{(a - c)^2 + (b - d)^2}{a^2 + b^2 + c^2 + d^2}}$.
Since $M^2$ is free of flat points $-1 <A< 1$, $A \neq 0$ for  $(u, v) \in {\mathcal D}_0$.
Hence, there exists  a function $\varphi$ such that $\tanh 4 \varphi = A$.
Then from \eqref{Eq-14} it follows that
$\overline{a}\,\overline{c} + \overline{b}\,\overline{d} = 0$,
i.e. $\langle \sigma(\overline{x},\overline{x}), \sigma(\overline{x},\overline{y}) \rangle=0$.

Finally  we obtain that at any point of a  timelike surface free of flat points with $H=0$
we can introduce  a special tangent frame field $\{x, y\}$ and a special normal frame field $\{n_1, n_2\}$,
such that
\begin{equation}\label{Eq-3-a}
\begin{array}{l}
\vspace{2mm}
\sigma(x,x) = \nu \,n_1; \\
\vspace{2mm}
\sigma(x,y) = \qquad \mu\,n_2;  \\
\vspace{2mm}
 \sigma(y,y) = \nu \,n_1,
\end{array}
\end{equation}
where $\nu = \langle \sigma(x,x), n_1 \rangle$,
$\mu = \langle \sigma(x,y), n_2 \rangle$. Since $M^2$ is free of flat points, $\nu  \mu \neq 0$.

The tangent directions determined by the tangent vector fields $x$ and $y$
we call \emph{canonical directions} of the surface.
The canonical directions  are uniquely determined at any point of a  timelike surface with $H=0$
free of flat points.
Obviously, the normal frame field $\{n_1, n_2\}$ is also uniquely determined by the canonical tangents,
and the functions $\nu$, $\mu$ are geometric functions (invariants) of the surface.
We call this orthonormal frame field
$\{x, y, n_1, n_2\}$ the \emph{geometric frame field} of $M^2$.

With respect to the geometric  frame field
$\{x, y, n_1, n_2\}$ we have the following  Frenet-type derivative
formulas of $M^2$:
\begin{equation} \label{Eq-5-a}
\begin{array}{ll}
\vspace{2mm} \nabla'_xx=\quad \quad \quad \gamma_1\,y+\,\nu\,n_1;
& \qquad
\nabla'_x n_1=  \nu\,x \quad \quad \quad \quad \quad \quad  + \beta_1\,n_2;\\
\vspace{2mm} \nabla'_xy=\gamma_1\,x\quad \quad \qquad  \qquad + \mu\,n_2;  & \qquad
\nabla'_y n_1=\quad \quad \; - \nu \,y \quad\quad \quad  + \beta_2\,n_2;\\
\vspace{2mm} \nabla'_yx=\quad\quad \, -\gamma_2\,y \quad \quad \quad   +\mu\,n_2;  & \qquad
\nabla'_xn_2= \quad \quad \, - \mu \,y  -\beta_1\,n_1;\\
\vspace{2mm}
\nabla'_yy=- \gamma_2\,x \quad\quad\;\;\,
+ \nu\,n_1; & \qquad
\nabla'_y n_2= \mu \, x  \quad \quad\;\;\;  \,-\beta_2\,n_1,
\end{array}
\end{equation}
where  $\gamma_1 = y(\ln \sqrt{-E}), \,\, \gamma_2 = - x(\ln
\sqrt{G})$, $\beta_1 =  \langle \nabla'_x n_1, n_2\rangle$, $\beta_2 =
\langle \nabla'_y n_1, n_2\rangle$.

\vskip 1mm
We shall use the following terminology: the integral lines of the canonical tangent
directions are said to be the \emph{canonical lines}; any parameters
$(u,v)$ of $M^2$ generating canonical parametric lines are said to be
\textit{semi-canonical parameters}.
It is clear that semi-canonical parameters are determined up to  changes of the following form:
$ u=u(\bar u), \,\, v=v(\bar v)$ or $ u=u(\bar v), \,\, v=v(\bar u)$.

The functions $\nu,\, \mu,\, \gamma_1,\, \gamma_2,\, \beta_1, \, \beta_2$ in the Frenet-type
formulas \eqref{Eq-5-a} are invariants of the surface $M^2$ determined by the canonical tangent directions.

\vskip 1mm
Using the Gauss and Codazzi equations from \eqref{Eq-5-a} we obtain  the following
equalities for the invariants of $M^2$:
\begin{equation} \label{Eq-6-a}
\begin{array}{l}
2\mu \, \gamma_2 + \nu \,\beta_2 =  x(\mu),\\
[2mm]
- 2 \mu \, \gamma_1 + \nu\,\beta_1 =  y(\mu),\\
[2mm]
2\nu \, \gamma_2 - \mu\,\beta_2 = x(\nu),\\
[2mm]
- 2\nu \, \gamma_1 -  \mu\,\beta_1 =  y(\nu),\\
[2mm]
\nu^2 - \mu^2 = x(\gamma_2) + y(\gamma_1) + \gamma_1^2 - \gamma_2^2,\\
[2mm]
 2 \nu \mu = x(\beta_2) - y(\beta_1) - \gamma_1 \,\beta_1 - \gamma_2 \,\beta_2.
\end{array}
\end{equation}

\vskip 2mm
The Gauss curvature of $M^2$ is $K = \langle \sigma(x,x), \sigma(y,y) \rangle - \langle \sigma(x,y), \sigma(x,y) \rangle$.
Hence, using \eqref{Eq-3-a} we get
\begin{equation} \label{Eq-7-a}
K = \nu^2 - \mu^2.
\end{equation}

\vskip 3mm
\begin{prop}\label{P:flat surface}
Let $M^2$ be a timelike surface with zero mean curvature  free of flat points.
Then at each point of the surface the Gauss curvature $K$ is non-zero.
\end{prop}
\noindent
\emph{Proof:}
Suppose that $K = 0$, i.e. $\mu^2 = \nu^2$.
From the first four equalities of \eqref{Eq-6-a} we get $\beta_1 = \beta_2 = 0$. Then the last equality of \eqref{Eq-6-a}
implies   $\nu \mu =0$, which contradicts the assumption that $M^2$ is free of flat points.
\qed

\vskip 3mm
The curvature of the normal connection (the normal curvature) of $M^2$ is given by the formula
$\varkappa= \langle R^{\bot}(x,y)n_2, n_1 \rangle$, where
$R^{\bot}(x,y)$ stands for the  normal curvature operator:
$$R^{\bot}(x,y) = D_xD_y- D_yD_x - D_{[x,y]}.$$
Using \eqref{Eq-5-a} and \eqref{Eq-6-a} we obtain that
\begin{equation}\label{Eq-8-a}
\varkappa = - 2 \nu \mu.
\end{equation}

Note that $\nu \mu \neq 0$ since $M^2$ is free of flat points. Thus the following statement holds true.

\begin{prop}\label{P:flat normal connection}
Let $M^2$ be a  timelike surface  with zero mean curvature  free of flat points. Then at each point of the surface
the curvature of the normal connection  $\varkappa$ is non-zero.
\end{prop}

Using that $\nu^2 \neq \mu^2$ from \eqref{Eq-6-a} it follows that $\gamma_1^2 + \gamma_2^2 \neq 0$, \,
$\beta_1^2 + \beta_2^2 \neq 0$.

\vskip 2mm
Now we shall introduce canonical parameters on each  timelike surface  with zero mean curvature  free of flat points.

\begin{defn}\label{D:Canonocal parameters-a}
Let $M^2$ be a timelike surface  with zero mean curvature  free of flat points and $(u,v)$ be
semi-canonical parameters. The parameters $(u,v)$ are said to be \textit{canonical}, if
$$E = \ds{- \frac{1}{\sqrt{\mu^2 + \nu^2}}};\qquad F = 0; \qquad G = \ds{\frac{1}{\sqrt{\mu^2 + \nu^2}}}.$$
\end{defn}

The canonical parameters are special isothermal parameters of the surface.

\begin{thm}\label{T:Canonocal parameters-a}
Any  timelike surface with zero mean curvature free of flat points locally admits canonical parameters.
\end{thm}

\noindent
\emph{Proof:}
Let $M^2$ be a timelike surface with $H=0$  free of flat points. We assume that $M^2$ is parameterized by
semi-canonical parameters $(u,v)$.
From equalities \eqref{Eq-6-a}, using that $z_u = \sqrt{-E}\,x, \,\, z_v= \sqrt{G}\, y$, we get
$$\gamma_1 = \ds{-\frac{1}{4 \sqrt{G}}\left(\ln (\mu^2 + \nu^2)\right)_v}; \qquad
\gamma_2 = \ds{\frac{1}{4 \sqrt{-E}}\left(\ln (\mu^2 + \nu^2)\right)_u}.$$
On the other hand, $\gamma_1 = \ds{y(\ln \sqrt{-E}) = \frac{1}{\sqrt{G}} \left(\ln\sqrt{-E}\right)_v}, \,\,
\gamma_2 =  \ds{ - x(\ln \sqrt{G}) = - \frac{1}{\sqrt{-E}} \left(\ln\sqrt{G}\right)_u}$.
Hence,
$$\ds{\frac{\partial}{\partial v}\left(\ln \left(E^2 (\mu^2 + \nu^2)\right)\right)} = 0;
\qquad \ds{\frac{\partial}{\partial u}\left(\ln \left(G^2 (\mu^2 + \nu^2)\right)\right)} = 0,$$
which imply that
$E \sqrt{\mu^2 + \nu^2}$ does not depend on $v$, and $G \sqrt{\mu^2 + \nu^2}$
does not depend on $u$.
Hence,
there exist functions $\varphi(u) >0 $ and $ \psi(v) >0$, such that
$$- E \sqrt{\mu^2 + \nu^2} = \varphi(u); \qquad G \sqrt{\mu^2 + \nu^2} = \psi(v).$$
Under the following change of the parameters:
$$\begin{array}{l}
\vspace{2mm}
\overline{u} = \ds{\int_{u_0}^u \sqrt{\varphi(u)}\, du}
+ \overline{u}_0, \quad \overline{u}_0 = const\\
[2mm]
\overline{v} = \ds{\int_{v_0}^v  \sqrt{\psi(v)}\, dv + \overline{v}_0},
\quad \overline{v}_0 = const
\end{array}$$
we obtain
$$\overline{E} = \ds{- \frac{1}{ \sqrt{\mu^2 + \nu^2}}}; \qquad \overline{F} = 0;
\qquad \overline{G} = \ds{\frac{1}{ \sqrt{\mu^2 + \nu^2}}},$$
i.e. the parameters $(\overline{u},\overline{v})$ are canonical. \qed

\vskip 3mm
Now let $M^2: z = z(u,v), \,\, (u,v) \in {\mathcal D}$ be a timelike surface  with $H=0$ free of flat points
and $(u,v)$ be canonical parameters, i.e.
$$E = \ds{- \frac{1}{ \sqrt{\mu^2 + \nu^2}}}; \qquad
F = 0; \qquad G = \ds{\frac{1}{ \sqrt{\mu^2 + \nu^2}}}.$$
Then the functions $\gamma_1$, $\gamma_2$, $\beta_1$ and $\beta_2$
are expressed as follows:
$$\gamma_1 = - \left((\mu^2 + \nu^2)^{\frac{1}{4}}\right)_v; \qquad \gamma_2
= \left((\mu^2 + \nu^2)^{\frac{1}{4}}\right)_u;$$
$$\beta_1 = \ds{(\mu^2 + \nu^2)^{\frac{1}{4}}
\left(\arctan \frac{\mu}{\nu}\right)_v}; \qquad \beta_2 =
\ds{(\mu^2 + \nu^2)^{\frac{1}{4}} \left(\arctan \frac{\mu}{\nu}\right)_u}.$$
The last two equalities of  \eqref{Eq-6-a} imply the following partial differential equations:
\begin{equation}\label{Eq-9-aa}
\begin{array}{l}
\ds{(\mu^2 + \nu^2)^{\frac{1}{2}}\, \Delta^h \ln (\mu^2 + \nu^2)^{\frac{1}{4}} =  \nu^2 - \mu^2};\\
[2mm]
\ds{(\mu^2 + \nu^2)^{\frac{1}{2}}\, \Delta^h \arctan \frac{\mu}{\nu} = 2 \nu  \mu}, \end{array}
\end{equation}
where $\Delta^h = \ds{\frac{\partial^2}{\partial u^2} - \frac{\partial^2}{\partial v^2}}$
is the hyperbolic Laplace operator.

\vskip 3mm
Now we shall prove the fundamental theorem of Bonnet-type for the class of
timelike surfaces with zero mean curvature free of flat points.

\begin{thm}\label{T:Fundamental Theorem - timelike}
Let $\mu(u,v)$ and $\nu(u,v)$ be two smooth functions, defined in a domain
${\mathcal D}, \,\, {\mathcal D} \subset {\R}^2$, and satisfying the conditions
$$\begin{array}{l}
\mu \nu \neq 0, \quad \mu^2-\nu^2\neq 0;\\
[2mm]
\ds{(\mu^2 + \nu^2)^{\frac{1}{2}}\, \Delta^h \ln (\mu^2 + \nu^2)^{\frac{1}{4}} =  \nu^2 - \mu^2};\\
[2mm]
\ds{(\mu^2 + \nu^2)^{\frac{1}{2}}\, \Delta^h \arctan \frac{\mu}{\nu} = 2 \nu  \mu}.
\end{array} $$
 Let $\{x_0, \, y_0, \, (n_1)_0,\, (n_2)_0\}$ be an
orthonormal frame at a point $p_0 \in \R^4_1$ such that $x_0$ is timelike and $y_0, \, (n_1)_0, (n_2)_0$ are spacelike.
Then there exist a subdomain ${\mathcal{D}}_0
\subset \mathcal{D}$ and a unique  timelike surface
$M^2: z = z(u,v), \,\, (u,v) \in {\mathcal D}_0, \,\, {\mathcal D}_0 \subset {\mathcal D}$  with zero mean curvature
free of flat points,
such that $M^2$ passes through $p_0$, the functions   $\nu(u,v)$, $\mu(u,v)$ are the geometric functions
of $M^2$ and $\{x_0, \, y_0, \, (n_1)_0,\, (n_2)_0\}$ is the geometric
frame of $M^2$ at the point $p_0$.
Furthermore, $(u,v)$ are canonical parameters of $M^2$.
\end{thm}

\vskip 2mm \noindent \emph{Proof:} Let us denote $E = -(\mu^2 + \nu^2)^{-\frac{1}{2}}; \,\,
G = (\mu^2 + \nu^2)^{-\frac{1}{2}}; \,\, \gamma_1 = -\left((\mu^2 +
\nu^2)^{\frac{1}{4}}\right)_v; \,\, \gamma_2 = \left((\mu^2 +
\nu^2)^{\frac{1}{4}}\right)_u$; $\beta_1 = \ds{(\mu^2 +
\nu^2)^{\frac{1}{4}} \left(\arctan \frac{\mu}{\nu}\right)_v};
\,\, \beta_2 =  \ds{(\mu^2 + \nu^2)^{\frac{1}{4}} \left(\arctan
\frac{\mu}{\nu}\right)_u}$.
 We consider the following
system of partial differential equations for the unknown vector
functions $x = x(u,v), \, y = y(u,v), \,n_1 = n_1(u,v), \,n_2 = n_2(u,v)$
in $\R^4_1$:
\begin{equation}\label{Eq-10}
\begin{array}{ll}
\vspace{2mm} x_u = \sqrt{-E}\, \gamma_1\, y + \sqrt{-E}\, \nu\, n_1
& \qquad x_v =
- \sqrt{G}\, \gamma_2\, y  + \sqrt{G}\, \mu\, n_2\\
\vspace{2mm} y_u =  \sqrt{-E}\, \gamma_1\, x + \sqrt{-E}\, \mu\, n_2  &
\qquad y_v = -\sqrt{G}\, \gamma_2\, x + \sqrt{G}\, \nu\, n_1 \\
\vspace{2mm} (n_1)_u = \sqrt{-E}\, \nu\, x + \sqrt{-E}\, \beta_1\, n_2  &
\qquad (n_1)_v = - \sqrt{G}\, \nu\, y + \sqrt{G}\, \beta_2\, n_2 \\
\vspace{2mm} (n_2)_u = - \sqrt{-E}\, \mu\, y - \sqrt{-E}\, \beta_1\, n_1 &
\qquad (n_2)_v = \sqrt{G}\, \mu\, x - \sqrt{G}\, \beta_2\, n_1
\end{array}
\end{equation}
We denote
$$Z =
\left(%
\begin{array}{c}
  x \\
  y \\
  n_1 \\
  n_2 \\
\end{array}%
\right); \quad
A = \sqrt{-E} \left(%
\begin{array}{cccc}
  0 & \gamma_1 & \nu & 0 \\
  \gamma_1 & 0 & 0 &  \mu \\
  \nu & 0 & 0 &  \beta_1 \\
  0 & -\mu & -\beta_1 & 0 \\
\end{array}%
\right); \quad B = \sqrt{G}
\left(%
\begin{array}{cccc}
  0 & -\gamma_2 & 0 &  \mu \\
  -\gamma_2 & 0 & \nu & 0 \\
  0 & -\nu & 0 &  \beta_2 \\
  \mu & 0 & -\beta_2 & 0 \\
\end{array}%
\right).$$
Then system \eqref{Eq-10} can be rewritten in the form:
\begin{equation}\label{Eq-11}
\begin{array}{l}
\vspace{2mm}
Z_u = A\,Z,\\
\vspace{2mm} Z_v = B\,Z.
\end{array}
\end{equation}
The integrability conditions of \eqref{Eq-11} are
$$Z_{uv} = Z_{vu},$$
i.e.
\begin{equation}\label{Eq-12}
\displaystyle{\frac{\partial a_i^k}{\partial v} - \frac{\partial b_i^k}{\partial u}
+ \sum_{j=1}^{4}(a_i^j\,b_j^k - b_i^j\,a_j^k) = 0, \quad i,k = 1,
\dots, 4,}
\end{equation}
 where $a_i^j$ and $b_i^j$ are the
elements of the matrices $A$ and $B$. Using the expressions of $E, G, \gamma_1, \gamma_2, \beta_1, \beta_2$
in terms of $\mu$, $\nu$ and equalities \eqref{Eq-9-aa}
we obtain that the integrability conditions
\eqref{Eq-12} are fulfilled. Hence, there exist a subset
$\mathcal{D}_1 \subset \mathcal{D}$ and unique vector functions $x
= x(u,v), \, y = y(u,v), \,n_1 = n_1(u,v), \,n_2 = n_2(u,v), \,\, (u,v)
\in \mathcal{D}_1$, which satisfy system \eqref{Eq-10} and the initial conditions
$$x(u_0,v_0) = x_0, \quad y(u_0,v_0) = y_0, \quad n_1(u_0,v_0) = (n_1)_0, \quad n_2(u_0,v_0) = (n_2)_0.$$

We shall prove that $x(u,v), \, y(u,v), \,n_1(u,v), \,n_2(u,v)$ form
an orthonormal frame in $\R^4_1$ for each $(u,v) \in
\mathcal{D}_1$. Let us consider the following functions:
$$\begin{array}{lll}
\vspace{2mm}
  \varphi_1 = \langle x,x \rangle + 1; & \qquad \varphi_5 =
  \langle x,y \rangle; & \qquad \varphi_8 = \langle y,n_1 \rangle; \\
\vspace{2mm}
  \varphi_2 = \langle y, y \rangle - 1; & \qquad \varphi_6 =
  \langle x,n_1 \rangle; & \qquad \varphi_9 = \langle y,n_2 \rangle; \\
\vspace{2mm}
  \varphi_3 = \langle n_1, n_1 \rangle - 1; & \qquad \varphi_7 =
  \langle x,n_2 \rangle; & \qquad \varphi_{10} = \langle n_1,n_2 \rangle; \\
\vspace{2mm}
  \varphi_4 = \langle n_2,n_2 \rangle - 1; &   &  \\
\end{array}$$
defined for each $(u,v) \in \mathcal{D}_1$. Using that $x(u,v), \,
y(u,v), \,n_1(u,v), \,n_2(u,v)$ satisfy \eqref{Eq-10}, we obtain  the system
\begin{equation}\label{Eq-12-a}
\begin{array}{lll}
\vspace{2mm}
\displaystyle{\frac{\partial \varphi_i}{\partial u} = \alpha_i^j \, \varphi_j},\\
\vspace{2mm} \displaystyle{\frac{\partial \varphi_i}{\partial v} =
\beta_i^j \, \varphi_j};
\end{array} \qquad i = 1, \dots, 10,
\end{equation}
where $\alpha_i^j, \beta_i^j, \,\, i,j = 1, \dots, 10$ are
functions of $(u,v) \in \mathcal{D}_1$. System \eqref{Eq-12-a} is a linear
system of partial differential equations for the functions
$\varphi_i(u,v), \,\,i = 1, \dots, 10, \,\,(u,v) \in
\mathcal{D}_1$, satisfying $\varphi_i(u_0,v_0) = 0, \,\,i = 1,
\dots, 10$. Hence, $\varphi_i(u,v) = 0, \,\,i = 1, \dots, 10$ for
each $(u,v) \in \mathcal{D}_1$. Consequently, the vector functions
$x(u,v), \, y(u,v), \,n_1(u,v), \,n_2(u,v)$ form an orthonormal frame
in $\R^4_1$ for each $(u,v) \in \mathcal{D}_1$.

Now, let us consider the  system
\begin{equation}\label{Eq-13}
\begin{array}{lll}
\vspace{2mm}
z_u = \sqrt{-E}\, x\\
\vspace{2mm} z_v = \sqrt{G}\, y
\end{array}
\end{equation}
of partial differential equations for the vector function
$z(u,v)$. Using that $x(u,v)$ and $y(u,v)$ satisfy  \eqref{Eq-10} we get that the integrability
conditions $z_{uv} = z_{vu}$ of system \eqref{Eq-13}
 are fulfilled. Hence,  there exist a subset $\mathcal{D}_0 \subset \mathcal{D}_1$ and
a unique vector function $z = z(u,v)$, defined for $(u,v) \in
\mathcal{D}_0$ and satisfying $z(u_0, v_0) = p_0$.

Consequently, the surface $M^2: z = z(u,v), \,\, (u,v) \in
\mathcal{D}_0$ satisfies the assertion of the theorem. \qed

\vskip 3mm
The meaning of Theorem \ref{T:Fundamental Theorem - timelike} is that
\textit{any timelike
 surface with $H=0$ free of flat points is determined up to a motion in $\R^4_1$ by two invariant functions
satisfying the system of natural partial differential equations  \eqref{Eq-9-aa}}.

\section{Rotational surfaces of Moore type with zero mean curvature} \label{S:Examples}
Considering general rotations in the Euclidean space $\R^4$, in \cite{M} C. Moore introduced
general rotational surfaces
as follows. Let $m: x(u) = \left( x^1(u), x^2(u),  x^3(u),
x^4(u)\right); \,\, u \in J \subset \R$ be a smooth curve in
$\R^4$, and $\alpha$, $\beta$ be constants. A general rotation of
the meridian curve  $m$ in $\R^4$ is given by
$$X(u,v)= \left( X^1(u,v), X^2(u,v),  X^3(u,v), X^4(u,v)\right),$$
where
$$\begin{array}{ll}
\vspace{2mm} X^1(u,v) = x^1(u)\cos\alpha v - x^2(u)\sin\alpha v; &
\qquad
X^3(u,v) = x^3(u)\cos\beta v - x^4(u)\sin\beta v; \\
\vspace{2mm} X^2(u,v) = x^1(u)\sin\alpha v + x^2(u)\cos\alpha v;&
\qquad X^4(u,v) = x^3(u)\sin\beta v + x^4(u)\cos\beta v.
\end{array}$$
In the case $\beta = 0$ the $X^3X^4$-plane is fixed and one gets
the classical rotation about a fixed two-dimensional axis.

In \cite{Mil}   we considered a surface $\mathcal{M}$ in  $\R^4$,
defined by the vector-valued function
$$z(u,v) = \left( f(u) \cos\alpha v, f(u) \sin \alpha v, g(u) \cos \beta v, g(u) \sin \beta v \right);
\quad u \in J \subset \R, \,\,  v \in [0; 2\pi),$$
where $f(u)$ and $g(u)$ are smooth functions, satisfying $\alpha^2
f^2(u)+ \beta^2 g^2(u) > 0$, $f'\,^2(u)+ g'\,^2(u) > 0$, $u \in
J,$ and $\alpha, \beta$ are positive constants, $\alpha \neq \beta$. The surface
$\mathcal{M}$  is a general rotational surface
whose meridians lie in  two-dimensional planes. In this case the
meridian is given by $m: x(u) = \left( f(u), 0,  g(u), 0\right);
\,\, u \in J \subset \R$. We found all minimal superconformal general
rotational surfaces in $\R^4$ \cite{Mil}.

\vskip 2mm
 Now let $m: x(u) = \left( x^1(u), x^2(u),  x^3(u),
x^4(u)\right); \,\, u \in J \subset \R$ be a spacelike or timelike
curve in the Minkowski space $\R^4_1$. We assume that $m$ is
parametrized by the arclength and denote $(\dot{x}^1)^2 +
(\dot{x}^2)^2 + (\dot{x}^3)^2 - (\dot{x}^4)^2 = \varepsilon$, where
$\varepsilon = 1$ if $c$ is spacelike and $\varepsilon = -1$ if
$c$ is timelike. Using the idea of Moore we consider the
surface $\mathcal{M}$ given by
$$\mathcal{M}: X(u,v)= \left( X^1(u,v), X^2(u,v),
X^3(u,v), X^4(u,v)\right),$$ where
$$\begin{array}{ll}
\vspace{2mm} X^1(u,v) = x^1(u)\cos\alpha v - x^2(u)\sin\alpha v; &
\qquad
X^3(u,v) = x^3(u)\cosh\beta v + x^4(u)\sinh\beta v; \\
\vspace{2mm} X^2(u,v) = x^1(u)\sin\alpha v + x^2(u)\cos\alpha v;&
\qquad X^4(u,v) = x^3(u)\sinh\beta v + x^4(u)\cosh\beta v,
\end{array}$$
$\alpha$, $\beta$ are  constants. This surface can be considered as rotational surface of
Moore type in  $\R^4_1$ obtained  by the meridian curve $m$.

Note that if $\beta = 0$, $\alpha = 1$ and $x^2(u) = 0$ (or $x^1(u)=0$),
one gets the surface
$$z(u,v) = \left(x^1(u) \cos v, x^1(u) \sin v, x^3(u), x^4(u)\right);
\quad u \in J,\,\,  v \in [0; 2\pi).$$
This is a  surface obtained by the rotation of the meridian curve about the two-dimensional
Lorentz plane $Oe_3e_4$. It is called a \emph{rotational surface of elliptic type}.
A local classification of spacelike surfaces in $\R^4_1$, which are invariant under spacelike rotations,
and with  mean curvature vector either vanishing or lightlike, is obtained in \cite{Haesen-Ort-2}.

If $\alpha = 0$, $\beta = 1$ and $x^3(u) = 0$ (or $x^4(u)=0$) one  gets \emph{a rotational surface of hyperbolic type}:
$$z(u,v) = \left( x^1(u), x^2(u), x^4(u) \sinh v, x^4(u) \cosh v\right);
\quad u \in J,\,\,  v \in \R.$$
It is an orbit of a smooth curve under the action of the orthogonal transformations of $\R^4_1$
which leave the spacelike plane $Oe_1e_2$ point-wise fixed.
A local classification of spacelike rotational surfaces of hyperbolic type
with vanishing or lightlike mean curvature vector field is given in \cite{Haesen-Ort-1}.
A classification of all timelike and spacelike hyperbolic rotational surfaces
with non-zero constant mean curvature in the three-dimensional de Sitter space
$\mathbb{S}^3_1$ is given in \cite{Liu-Liu-1}.
Similarly, a classification of the spacelike and timelike Weingarten rotational surfaces in $\mathbb{S}^3_1$
is found in \cite{Liu-Liu-2}.

\vskip 2mm
In \cite{GM5} we considered a special class of spacelike
surfaces  in the Minkowski space $\R^4_1$ parameterized by
\begin{equation} \label{Eq-example}\notag
\mathcal{M}_1: z(u,v) = \left( f(u) \cos \alpha v, f(u) \sin \alpha v, g(u) \cosh \beta v, g(u) \sinh \beta v \right),
\end{equation}
where $f(u)$ and $g(u)$ are smooth functions, satisfying
$$\alpha^2 f^2(u)- \beta^2 g^2(u) > 0, \quad f'\,^2(u)+ g'\,^2(u) > 0, \quad u \in J
\subset \R,$$
$v \in [0; 2\pi)$, and $\alpha, \beta$ are positive constants, $\alpha \neq \beta$.
These surfaces are rotational surfaces of Moore type obtained by
 the meridian curve $m: x(u) = \left( f(u), 0,  g(u), 0\right)$ lying in the two-dimensional plane $Oe_1e_3$.

\vskip 2mm
Now we shall consider  rotational surfaces of Moore type in $\R^4_1$ obtained by the
meridian curve $m: x(u) = \left( f(u), 0,  0, g(u)\right)$ lying in the Lorentz  plane $Oe_1e_4$.
The rotational surface is spacelike if the meridian $m$ is spacelike, and the rotational surface is timelike if
 $m$ is timelike.
 We shall consider the timelike case and we shall find all timelike rotational surface of Moore type with zero mean curvature.

Let $\mathcal{M}_2$ be the surface parameterized by
\begin{equation} \label{Eq-example}
\mathcal{M}_2: z(u,v) = \left( f(u) \cos \alpha v, f(u) \sin \alpha v, g(u) \sinh \beta v, g(u) \cosh \beta v \right),
\end{equation}
where $f(u)$ and $g(u)$ are smooth functions, satisfying
$$f'\,^2(u)- g'\,^2(u) < 0, \quad \alpha^2 f^2(u)+ \beta^2 g^2(u) > 0,  \quad u \in J
\subset \R,$$
$v \in [0; 2\pi)$, and $\alpha, \beta$ are positive constants, $\alpha \neq \beta$.

The coefficients of the first fundamental form of $\mathcal{M}_2$
are
\begin{equation} \label{Eq-example-1}
E = f'\,^2(u)- g'\,^2(u); \qquad F = 0; \qquad G =\alpha^2 f^2(u)+ \beta^2
g^2(u).
\end{equation}
Hence, $\mathcal{M}_2$ is a timelike surface.

We consider
the following normal frame field of $\mathcal{M}_2$:
$$\begin{array}{l}
\vspace{2mm}
n_1 = \ds{\frac{1}{\sqrt{g'\,^2 - f'\,^2}}\left(g' \cos \alpha v, g' \sin \alpha v, f' \sinh \beta v, f' \cosh \beta v \right)};\\
\vspace{2mm} n_2 = \ds{\frac{1}{\sqrt{\alpha^2 f^2 + \beta^2
g^2}}\left(  \beta g \sin \alpha v, - \beta g \cos \alpha v, \alpha
f \cosh \beta v,  \alpha f \sinh \beta v \right)}.
\end{array}$$
We have that $\langle n_1, n_1 \rangle = 1,\,\, \langle n_2, n_2
\rangle =1, \,\, \langle n_1, n_2 \rangle =0$. Denoting $x = \ds{\frac{z_u}{\sqrt{g'\,^2- f'\,^2}}}$,
$y = \ds{\frac{z_v} {\sqrt{\alpha^2 f^2+ \beta^2 g^2}}}$, we obtain geometric  frame field
$\{x, y, n_1, n_2\}$.

The mean curvature vector field  $H$ of $\mathcal{M}_2$ is given by  the following formula:
$$H = \ds{\frac{1}{2} \left(- \frac{g' f'' - f' g''}{(g'\,^2 -
f'\,^2)^{\frac{3}{2}}} -  \frac{\alpha^2 f g' + \beta^2 g f'}{\sqrt{g'\,^2 - f'\,^2}(\alpha^2 f^2 + \beta^2 g^2)} \right) n_1}.$$
Hence, $\mathcal{M}_2$ has zero mean curvature  ($H = 0$) if and only if  the functions $f(u)$ and $g(u)$ satisfy the equality
\begin{equation}\label{Eq-minimal}
\ds{- \frac{g' f'' - f' g''}{g'\,^2 -
f'\,^2} =  \frac{\alpha^2 f g' + \beta^2 g f'}{\alpha^2 f^2+ \beta^2 g^2}}.
\end{equation}

We shall find all timelike  rotational surfaces of Moore type with $H=0$.

\vskip 2mm
\begin{prop}
The timelike surface $\mathcal{M}_2$ defined by \eqref{Eq-example}  has zero mean curvature if and only if the functions
$f(u)$ and $g(u)$ satisfy
$$f'\,^2 = \frac{A - \alpha^2 f^2}{\alpha^2f^2 + \beta^2 g^2}; \qquad g'\,^2 = \frac{A + \beta^2 g^2}{\alpha^2f^2 + \beta^2 g^2},$$
where $A$ is a positive constant.
\end{prop}

\vskip 2mm \noindent \emph{Proof:}
Let $\mathcal{M}_2$ has zero mean curvature.
Calculating the functions $\nu$ and $\mu$  in the Frenet-type derivative
formulas \eqref{Eq-5-a} for this surface, we obtain:
\begin{equation} \label{Eq-example-2}
\begin{array}{ll}
\vspace{2mm} \nu = \ds{- \frac{\alpha^2 f g' + \beta^2 g f'}{\sqrt{g'\,^2 - f'\,^2}(\alpha^2 f^2 + \beta^2 g^2)}};& \qquad
 \mu = \ds{\frac{\alpha \beta (g' f - f' g)}{\sqrt{g'\,^2 - f'\,^2}(\alpha^2 f^2 + \beta^2 g^2)}}.
 \end{array}
\end{equation}

From the first and the third equalities of  \eqref{Eq-6-a} it follows that
$ \gamma_2 = \ds{\frac{1}{4}} \,x\left( \ln(\mu^2 + \nu^2)\right)$.
On the other hand we have $\gamma_2 = - x\left( \ln \sqrt{G}\right)$.
Hence, we get
$\ds{\frac{1}{4}} \,x\left( \ln(\mu^2 + \nu^2)\right) + \ds{\frac{1}{2}}\, x\left( \ln G\right)=  0$,
which implies that $$x\left( (\mu^2 + \nu^2) G^2\right) = 0.$$
Now, using that
$\mu$, $\nu$, and $G$ are functions depending only on the parameter $u$, we
obtain
\begin{equation} \label{Eq-example-3}
(\mu^2 + \nu^2) G^2 = c^2,
\end{equation}
where $c$ is a constant.
From \eqref{Eq-example-1}, \eqref{Eq-example-2}, and \eqref{Eq-example-3} it follows that
\begin{equation} \label{Eq-example-4}
\frac{\alpha^2 \beta^2 (g' f - f' g)^2 + (\alpha^2 f g' + \beta^2 g f')^2}{g'\,^2 - f'\,^2} = c^2.
\end{equation}

Without loss of generality we can assume that $g'\,^2 - f'\,^2 = 1$. Then \eqref{Eq-example-4} implies
\begin{equation} \label{Eq-example-5}
\alpha^2 f^2 g'\,^2 + \beta^2 g^2 f'\,^2 = \frac{c^2}{\alpha^2 + \beta^2}.
\end{equation}
We denote $A = \ds{\frac{c^2}{\alpha^2 + \beta^2}}$.
Now, using that $g'\,^2 = 1 + f'\,^2$, from \eqref{Eq-example-5} we obtain
\begin{equation} \label{Eq-example-6}
f'\,^2 = \frac{A - \alpha^2 f^2}{\alpha^2f^2 + \beta^2 g^2}; \qquad g'\,^2 = \frac{A +\beta^2 g^2}{\alpha^2f^2 + \beta^2 g^2}.
\end{equation}
If we assume that $A = 0$, then we get $f'\,^2 <0$ which is a contradiction. So, $A >0$.

Conversely, if $f(u)$ and $g(u)$ are functions satisfying \eqref{Eq-example-6}, then by a direct computation
it can be checked that
\eqref{Eq-example-5} is fulfilled and equality \eqref{Eq-minimal} holds true.
Hence, $\mathcal{M}_2$ is a  timelike surface with zero mean curvature.
\qed

\vskip 3mm
Note that the constant $A$ satisfies the inequality $ A > \alpha^2 f^2$, since $f'\,^2 >0$.

\vskip 2mm
Equalities \eqref{Eq-example-6} imply $(A - \alpha^2 f^2)g'\,^2 = (A + \beta^2 g^2) f'\,^2$, i.e.
$$\frac{f'}{\sqrt{A - \alpha^2 f^2 }} = \varepsilon \frac{g'}{\sqrt{A + \beta^2 g^2}}, \qquad \varepsilon =\pm 1.$$
Integrating the last equality we obtain
$$\int \frac{df}{\sqrt{A - \alpha^2 f^2}} = \varepsilon \int \frac{dg}{\sqrt{A + \beta^2 g^2}}.$$
Calculating the integrals we get
$$\arcsin \frac{\alpha f}{\sqrt{A}} =\varepsilon \frac{\alpha}{\beta} \ln \left| \beta g + \sqrt{\beta^2 g^2 + A} \right| + C, \qquad C = const.$$

In such a way we find all timelike rotational surfaces of Moore type with $H=0$. They are described in the following

\begin{thm}\label{T:rotational-Itype}
The timelike  rotational surface of Moore type  has zero mean curvature if and only if
the meridian curve is given by the formula
$$f = \frac{\sqrt{A}}{\alpha} \sin \left( \frac{\varepsilon \alpha}{\beta} \ln \left| \beta g + \sqrt{\beta^2 g^2 +A} \right| + C \right), \qquad C = const.$$
\end{thm}

\vskip 3mm
The timelike  rotational surface of Moore type  with $H=0$  gives us a solution to system \eqref{Eq-9-aa}.
Indeed,
if we denote $\varphi = \ds{ \frac{\varepsilon \alpha}{\beta} \ln \left| \beta g + \sqrt{\beta^2 g^2 + A} \right| + C}$,
then the meridian curve of the corresponding rotational surface
$\mathcal{M}_2^{0}$ is
given by $f =\ds{ \frac{\sqrt{A}}{\alpha}  \sin \varphi}$, where $A = \ds{\frac{c^2}{\alpha^2 + \beta^2}}$.
In such case  the coefficients $E$ and $G$ of the first fundamental form of $\mathcal{M}_2^{0}$  are $E = - 1$,
$G =A \sin^2 \varphi + \beta^2 g^2$. The parameters $(u,v)$ are not canonical, since
$$\mu^2 + \nu^2 = \ds{\frac{c^2}{(A \sin^2 \varphi + \beta^2 g^2)^2} = \frac{c^2}{G^2}}.$$
In terms of the parameters $(u,v)$ the functions $\mu$ and $\nu$ of $\mathcal{M}_2^{0}$  are expressed as follows:
$$\mu = \mu(u) = \ds{\frac{\alpha \beta (g' f - f' g)}{A \sin^2 \varphi + \beta^2 g^2}}; \qquad
\nu  = \nu(u) = \ds{- \frac{\alpha^2 f g' + \beta^2 g f'}{A \sin^2 \varphi + \beta^2 g^2}}.$$

Let us consider the following change of the parameters:
$$\begin{array}{l}
\vspace{2mm}
\overline{u} = \sqrt{c}\ds{\int_{0}^u \frac{du}{\sqrt{A \sin^2 \varphi + \beta^2 g^2}}};
\\
[2mm]
\overline{v} = \sqrt{c}\,v.
\end{array}$$
Then
$$\overline{E} = \ds{- \frac{1}{ \sqrt{\mu^2 + \nu^2}}}; \qquad \overline{F} = 0;
\qquad \overline{G} = \ds{\frac{1}{ \sqrt{\mu^2 + \nu^2}}},$$
i.e.  $(\overline{u},\overline{v})$ are canonical parameters of the surface $\mathcal{M}_2^{0}$.
Consequently, the functions $\mu = \mu(u(\overline{u}))$ and
 $\nu = \nu(u(\overline{u}))$ are solutions to system \eqref{Eq-9-aa}.

\section{Conclusion} \label{S:Conclusion}

In Section \ref{S:timelike} we proved that any  timelike surface with $H=0$ free of flat points
is determined up to a motion in $\R^4_1$ by two invariant functions   $\mu$ and $\nu$
satisfying a system of two partial differential equations, namely system \eqref{Eq-9-aa}.

Using equalities \eqref{Eq-7-a} and \eqref{Eq-8-a}  we can express the geometric functions $\mu$ and $\nu$ of a minimal timelike surface in terms of the Gauss curvature $K$
and the curvature of the normal connection $\varkappa$ as follows:
$$\mu^2 +\nu^2 = \sqrt{K^2 + \varkappa^2}; \qquad
\ds{\frac{\mu}{\nu} = \frac{K \pm \sqrt{K^2 + \varkappa^2}}{\varkappa}}.$$

So  equalities \eqref{Eq-9-aa}  can be rewritten as
\begin{equation}\label{Eq-9-b}
\begin{array}{l}
\ds{(K^2 + \varkappa^2)^{\frac{1}{4}}\, \Delta^h \ln (K^2 + \varkappa^2)^{\frac{1}{8}}} = K ;\\
[2mm]
\ds{(K^2 + \varkappa^2)^{\frac{1}{4}}\, \Delta^h \arctan \frac{\varkappa}{K} = 2 \varkappa}.
\end{array}
\end{equation}
Then the fundamental theorem  for  timelike surfaces with zero mean curvature
can be stated in terms of the Gauss curvature $K$
and the curvature of the normal connection $\varkappa$ as follows:

\begin{thm}\label{T:Fundamental Theorem 2 - timelike}
Let $K(u,v) \neq 0$ and $\varkappa(u,v) \neq 0$ be two smooth functions, defined in a domain
${\mathcal D}, \,\, {\mathcal D} \subset {\R}^2$, and satisfying the equalities
\begin{equation}\notag
\begin{array}{l}
\ds{(K^2 + \varkappa^2)^{\frac{1}{4}}\, \Delta^h \ln (K^2 + \varkappa^2)^{\frac{1}{8}}} = K ;\\
[2mm]
\ds{(K^2 + \varkappa^2)^{\frac{1}{4}}\, \Delta^h \arctan \frac{\varkappa}{K} = 2 \varkappa}.
\end{array}
\end{equation}
Then there exists a unique (up to a motion in $\R^4_1$)  timelike surface  with zero mean curvature
free of flat points  such that $K(u,v)$ and  $\varkappa(u,v)$ are the Gauss curvature and the curvature of the normal connection,
respectively, and $(u,v)$ are the canonical parameters.
\end{thm}

\vskip 3mm
If we set $$K = e^{2X} \cos Y; \quad \varkappa = e^{2X} \sin Y,$$
where $X = X(u,v)$, $Y=Y(u,v)$, then system \eqref{Eq-9-b} takes the form

\begin{equation}\label{Eq-9-tim}
\begin{array}{l}
\Delta^h X = 2 e^{X} \cos Y ;\\
[2mm]
\Delta^h Y = 2 e^{X} \sin Y .
\end{array}
\end{equation}

The study of timelike surfaces with $H=0$  free of flat points is equivalent to the study of the
solutions of system  \eqref{Eq-9-tim}.

In Section \ref{S:Examples} we proved that the geometric  functions $\mu = \mu(u(\overline{u}))$ and
 $\nu = \nu(u(\overline{u}))$ of the timelike  rotational surface $\mathcal{M}_2^0$
 give a solution $(\mu, \nu)$ to system \eqref{Eq-9-aa}.
Hence, the functions $K(\overline{u}) = \nu^2(\overline{u}) - \mu^2(\overline{u})$  and
$\varkappa(\overline{u}) = -2 \nu(\overline{u}) \mu(\overline{u})$ give a solution $(K, \varkappa)$ to  system  \eqref{Eq-9-b}.
Then we can find a solution $(X,Y)$ to system   \eqref{Eq-9-tim}.

\vskip 2mm

In \cite{Alias-Palmer} it is proved that the spacelike surfaces with $H=0$ in Minkowski space $\R^4_1$ are described
by a system similar to \eqref{Eq-9-b}:
\begin{equation}\label{Eq-9-a}
\begin{array}{l}
\ds{(K^2 + \varkappa^2)^{\frac{1}{4}}\, \Delta \ln (K^2 + \varkappa^2)^{\frac{1}{8}}} = K ;\\
[2mm]
\ds{(K^2 + \varkappa^2)^{\frac{1}{4}}\, \Delta \arctan \frac{\varkappa}{K} = 2 \varkappa},
\end{array}
\end{equation}
where $K$ and $\varkappa$ are the Gauss curvature and the normal curvature, respectively.

Substituting  $K = e^{2X} \cos Y; \,\, \varkappa = e^{2X} \sin Y$,
system \eqref{Eq-9-a} takes the form

\begin{equation}\label{Eq-9-sp}
\begin{array}{l}
\Delta X = 2 e^{X} \cos Y ;\\
[2mm]
\Delta Y = 2 e^{X} \sin Y .
\end{array}
\end{equation}
Thus the study of  spacelike surfaces with zero mean curvature free of flat points
is equivalent to the study of solutions of system \eqref{Eq-9-sp}.

The spacelike  rotational surfaces of Moore type  with two-dimensional meridians generate
a one-parameter family of solutions to system  \eqref{Eq-9-sp}  similarly to the considerations  in  Section \ref{S:Examples}.

\vskip 2mm System \eqref{Eq-0} obtained by de Azevedo Tribuzy and
Guadalupe  \cite{Guad-Trib} for minimal non-superconformal
surfaces in the Euclidean space $\R^4$
 can be rewritten in the form:

\begin{equation}\label{Eq-9-bb}
\begin{array}{l}
\ds{(K^2 - \varkappa^2)^{\frac{1}{4}}\, \Delta \ln (K^2 - \varkappa^2)^{\frac{1}{8}}} = K ;\\
[2mm]
\ds{(K^2 - \varkappa^2)^{\frac{1}{4}}\, \Delta \ln \frac{K - \varkappa}{K + \varkappa} = - 4\varkappa}.
\end{array}
\end{equation}

Substituting  $K = e^{2X} \cosh Y; \,\, \varkappa = e^{2X} \sinh Y$,
system \eqref{Eq-9-bb} takes the form
\begin{equation}\label{Eq-9-real}
\begin{array}{l}
\Delta X = 2 e^{X} \cosh Y ;\\
[2mm]
\Delta Y = 2 e^{X} \sinh Y .
\end{array}
\end{equation}

The general rotational surfaces of Moore in Euclidean space $\R^4$ with meridians lying in two-dimensional planes
 generate a one-parameter family of solutions  to system  \eqref{Eq-9-real}.

\vskip 2mm
Hence, the background systems of partial differential equations for surfaces with zero mean curvature in $\R^4$ and $\R^4_1$
are systems \eqref{Eq-9-real}, \eqref{Eq-9-sp}, and \eqref{Eq-9-tim}.

 \vskip 5mm \textbf{Acknowledgements:} The second author is
partially supported by "L. Karavelov" Civil Engineering Higher
School, Sofia, Bulgaria under Contract No 10/2011.

\end{document}